\documentclass[a4paper,11pt]{article}
\usepackage[latin1]{inputenc}
\usepackage[T1]{fontenc}
\usepackage[francais]{babel}
\usepackage{amssymb}
\usepackage{amsfonts, amsmath}
 \newtheorem{thm}{THÉORÈME}

 \newtheorem{defn}{DÉFINITION}
 \newtheorem{rem}{REMARQUE}
 
\title{\textbf{ÉTUDE DES SOLUTIONS MÉROMORPHES D'ÉQUATIONS DIFFÉRENTIELLES}}
\author{\textbf{A. Lesfari}
\\Département de Mathématiques
\\Faculté des Sciences
\\Université Chouaïb Doukkali
\\B.P. 20, El-Jadida, Maroc.
\\E-mail : Lesfariahmed@yahoo.fr , lesfari@ucd.ac.ma }
\date{}
\begin{document}
\maketitle Abstract. In this paper we shall study differential
equations in the complex domain. The method of indeterminate
coefficients and the majorant method lead to a proof of the
existence and uniqueness of meromorphic solution of differential
equations. We discuss their connection with the concept of
algebraic integrability systems.

\section{Position du probl\`{e}me}
Dans ce travail, nous envisagerons l'\'{e}tude des \'{e}quations
diff\'{e}rentielles dans le domaine complexe. Soit le syst\`{e}me
d'\'{e}quations diff\'{e}rentielles non-lin\'{e}aires
\begin{eqnarray}
\frac {dw_{1}}{dz}&=&f_{1}\left( z,w_{1},...,w_{n}\right)
,\nonumber\\
&\vdots&\\
\frac {dw_{n}}{dz}&=&f_{n}\left(
z,w_{1},...,w_{n}\right),\nonumber
\end{eqnarray}
o\`{u} $f_{1},...,f_{n}$ sont des fonctions de $n+1$ variables
complexes $z,w_{1},...,w_{n}$ et qui appliquent un domaine de
$\mathbb{C}^{n+1}$ dans $\mathbb{C}$. Le probl\`{e}me de Cauchy
consiste en la recherche d'une solution $\left( w_{1}\left(
z\right) ,...,w_{n}\left( z\right) \right) $ dans un voisinage
d'un point $z_{0},$ passant par le point donn\'{e} $\left(
z_{0},w_{1}^{0},...,w_{n}^{0}\right) $ c'est-\`{a}-dire
satisfaisant aux conditions initiales
\begin{eqnarray}
w_{1}(z_{0})&=&w_{1}^{0},\nonumber\\
&\vdots&\nonumber\\
w_{n}(z_{0})&=&w_{n}^{0}.\nonumber
\end{eqnarray}
Notons que le syst\`{e}me $\left( 1\right) $ peut s'\'{e}crire
sous forme vectorielle dans $\mathbb{C}^{n}$
$$\frac{dw}{dz}=f(z,w(z)) ,$$
en posant $w=$ $(w_{1},...,w_{n})$ et $f=(f_{1},...,f_{n}).$ Dans
ce cas, le probl\`{e}me de Cauchy consistera \`{a} d\'{e}terminer
la solution $w(z)$ telle que $$w(z_{0})
=w_{0}=(w_{1}^{0},...,w_{n}^{0}).$$

Commençons tout d'abord par décrire quelques résultats connus. On
sait que lorsque les fonctions $f_{1},...,f_{n}$ sont holomorphes
au voisinage du point $\left( z_{0},w_{1}^{0},...,w_{n}^{0}\right)
$ alors le probl\`{e}me de Cauchy admet une solution holomorphe et
une seule. Une question naturelle se pose: le probl\`{e}me de
Cauchy peut-il admettre quelque solution non holomorphe au
voisinage du point $\left( z_{0},w_{1}^{0},...,w_{n}^{0}\right) ?$
Lorsque les fonctions $f_{1},...,f_{n}$ sont holomorphes, la
r\'{e}ponse est n\'{e}gative. D'autres circonstances peuvent se
produire pour le problème de Cauchy relatif au système d'équations
différentielles $\left( 1\right) ,$ lorsque l'hypothèse
d'holomorphie relative aux fonctions $f_{1},...,f_{n}$ n'est plus
satisfaite au voisinage d'un point. On constate dans une telle
éventualité que les comportements des solutions peuvent revêtir
les aspects les plus divers. En général, les singularités des
solutions sont de deux types: mobiles ou fixes, suivant qu'elles
dépendent ou non des conditions initiales. Des résultats
importants ont été obtenus par Painlevé $[14]$. Supposons par
exemple que le système $(1)$ s'écrit sous la forme
\begin{eqnarray}
\frac{dw_{1}}{dz}&=&\frac{P_{1}(z,w_{1}w_{n})}{Q_{1}(z,w_{1},...,w_{n})},\nonumber\\
&\vdots&\nonumber\\
\frac{dw_{n}}{dz}&=&\frac{P_{n}(
z,w_{1},...,w_{n})}{Q_{n}(z,w_{1},...,w_{n})},\nonumber
\end{eqnarray}
avec $$P_{k}\left( z,w_{1},...,w_{n}\right) =\sum_{0\leq
i_{1},\ldots ,i_{n}\leq p}A_{i_{1},\ldots ,i_{n}}^{\left( k\right)
}\left( z\right) w_{1}^{i_{1}}...w_{n}^{i_{n}},\text{ }1\leq k\leq
n,$$$$Q_{k}\left( z,w_{1},...,w_{n}\right) =\sum_{0\leq
j_{1},\ldots ,j_{n}\leq q}B_{j_{1},\ldots ,j_{n}}^{\left( k\right)
}\left( z\right) w_{1}^{j_{1}}...w_{n}^{j_{n}},\text{ }1\leq k\leq
n,$$ des polynômes à plusieurs indéterminées $w_{1},...,w_{n}$ et
à coefficients algébriques en $z.$ On
sait\\
\null\hskip 0.7cm $\left( i\right) $ que les singularit\'{e}s
fixes sont constitu\'{e}es par quatre ensembles de points. Le
premier est l'ensemble des points singuliers des coefficients
$A_{i_{1},\ldots ,i_{n}}^{\left( k\right) }\left( z\right) ,$
$B_{j_{1},\ldots ,j_{n}}^{\left( k\right) }\left( z\right) $
intervenant dans les polyn\^{o}mes $P_{k}\left(
z,w_{1},...,w_{n}\right) $ et $Q_{k}\left(
z,w_{1},...,w_{n}\right) .$ En g\'{e}n\'{e}ral cet ensemble
contient le point $z=\infty .$ Le second ensemble est
constitu\'{e} des points $\alpha _{j}$ tels que $$Q_{k}\left(
z,w_{1},...,w_{n}\right) =0,$$ circonstance qui se produit si les
coefficients $B_{j_{1},\ldots ,j_{n}}^{\left( k\right) }\left(
z\right) $ s'annulent tous pour $z=\alpha _{j}.$ Le troisi\`{e}me
est l'ensemble des points $\beta _{l}$ tels que pour certaines
valeurs $\left( w_{1^{\prime }},...,w_{n^{\prime }}\right) $ de
$\left( w_{1},...,w_{n}\right) ,$ on ait $$P_{k}\left( \beta
_{l},w_{1^{\prime }},...,w_{n^{\prime }}\right) =Q_{k}\left( \beta
_{l},w_{1^{\prime }},...,w_{n^{\prime }}\right) =0.$$ D\`{e}s lors
les seconds membres du syst\`{e}me ci dessus se pr\'{e}sentent
sous la forme ind\'{e}termin\'{e}e $\frac{0}{0}$ aux points
$\left( \beta _{l},w_{1^{\prime }},...,w_{n^{\prime }}\right) .$
Enfin, l'ensemble des points $\gamma _{m}$ tels qu'il existe des
valeurs $u_{1},...,u_{n},$ pour lesquelles $$R_{k}\left( \gamma
_{m},u_{1},...,u_{n}\right) =S_{k}\left( \gamma
_{m},u_{1},...,u_{n}\right) =0,$$ o\`{u} $R_{k}$ et $S_{k}$ sont
des polyn\^{o}mes en $u_{1},...,u_{n}$ obtenus \`{a} partir de
$P_{k}$ et $Q_{k}$ en posant $$w_{1}=\frac{1}{u_{1}},\ldots
,w_{n}=\frac{1}{u_{n}}.$$ Chacun de ces ensembles ne comporte
qu'un nombre fini d'\'{e}l\'{e}ments. Les singularit\'{e}s fixes
du
syst\`{e}me en question sont en nombre fini.\\
\null\hskip 0.7cm $\left( ii\right) $ que les singularit\'{e}s
mobiles de solutions de ce syst\`{e}me sont des singularit\'{e}s
mobiles alg\'{e}briques: p\^{o}les et (ou) points critiques
alg\'{e}briques. Il n'y a pas de points singuliers essentiels pour
la solution $\left( w_{1},...,w_{n}\right) .$\\

Consid\'{e}rant le syst\`{e}me d'\'{e}quations diff\'{e}rentielles
$\left( 1\right) ,$ peut-on trouver des conditions suffisantes
d'existence et d'unicit\'{e} de solutions m\'{e}romorphes? Nous
\'{e}tablirons un th\'{e}or\`{e}me d'existence et d'unicit\'{e}
pour la solution du probl\`{e}me de Cauchy relatif au syst\`{e}me
d'\'{e}quations diff\'{e}rentielles $\left( 1\right) ,$ en faisant
appel \`{a} la m\'{e}thode des coefficients ind\'{e}termin\'{e}s.
La solution sera explicit\'{e}e sous la forme d'une s\'{e}rie de
Laurent. Il se posera d\`{e}s lors le probl\`{e}me de la
convergence. Celui-ci sera r\'{e}solu par la m\'{e}thode des
fonctions majorantes (pour cette notion voir par exemple [3, 7,
9]). Nombreux sont les probl\`{e}mes, aussi bien th\'{e}orique que
pratique, ou apparaissent des \'{e}quations diff\'{e}rentielles
dont le second membre n'est pas holomorphe. Nous verrons, dans la
derni\`{e}re section, que les solutions m\'{e}romorphes
d\'{e}pendant d'un nombre suffisant de param\`{e}tres libres
jouent un r\^{o}le crucial dans l'\'{e}tude des \'{e}quations
diff\'{e}rentielles dites alg\'{e}briquement int\'{e}grables.

\section{Existence et unicit\'{e} de solutions m\'{e}romorphes}

Dans ce qui suivra, nous envisagerons le probl\`{e}me de Cauchy
relatif au syst\`{e}me normal $\left( 1\right) $ dans
l'hypoth\`{e}se o\`{u} $f_{1},...,f_{n}$ ne dépendent pas
explicitement de $z$ c'est-à-dire
\begin{eqnarray}
\frac {dw_{1}}{dz}&=&f_{1}\left(w_{1},...,w_{n}\right)
,\nonumber\\
&\vdots&\\
\frac {dw_{n}}{dz}&=&f_{n}\left(w_{1},...,w_{n}\right),\nonumber
\end{eqnarray}
On suppose que $f_{1},...,f_{n}$ sont des fonctions rationnelles
en $w_{1},...,w_{n}$ et que le syst\`{e}me $\left( 2\right) $ est
quasi-homog\`{e}ne, c'est-\`{a}-dire ils existent des entiers
positifs $s_{1},...,s_{n}$ telles que
$$
f_{i}(\alpha ^{s_{1}}w_{1},...,\alpha ^{s_{n}}w_{n})=\alpha
^{s_{i}+1}f_{i}(w_{1},...,w_{n}) ,\quad 1\leq i\leq n,
$$
pour chaque constante non nulle $\alpha $. Autrement dit, le
syst\`{e}me $(2) $ est invariant par la transformation
$$z\longrightarrow \alpha ^{-1}z,\text{ }w_{1}\longrightarrow\alpha
^{s_{1}}w_{1},\ldots ,\text{ }w_{n}\longrightarrow \alpha
^{s_{n}}w_{n}.$$ Notons que si le d\'{e}terminant
\begin{equation}\label{eqn:euler}
\Delta \equiv \det (w_{j}\frac{\partial f_{i}}{\partial
w_{j}}-\delta _{ij}f_{i}) _{1\leq i,j\leq n},
\end{equation}
est non identiquement nul, alors le choix des nombres
$s_{1},...,s_{n}$ est unique.\\
Dans tout ce qui va suivre, nous supposerons, pour simplifier les
notations que $z_0=w_0=0,$ ce qui n'affecte pas la généralité des
résultats.

\begin{thm}
Supposons que
\begin{equation}\label{eqn:euler}
w_{i}=\frac{1}{z^{s_{i}}}\sum_{k=0}^{\infty }c_{i}^{(k)
}z^{k},\quad 1\leq i\leq n,
\end{equation}
où $c^{(0)}\neq 0$, soit la solution formelle en séries de
Laurent, obtenue par la méthode des coefficients indéteminés, du
syst\`{e}me quasi-homogène $(2)$. Alors, les coefficients
$c_{i}^{(0)}$ satisfont aux \'{e}quations non-lin\'{e}aires
\begin{equation}\label{eqn:euler}
s_{i}c_{i}^{(0)}+f_{i}(c_{1}^{(0)},...,c_{n}^{(0)})=0,
\end{equation}
où $1\leq i\leq n,$ tandis que $c_{i}^{(1)},c_{i}^{(2)},...$
satisfont chacun \`{a} un syst\`{e}me d'\'{e}quations
lin\'{e}aires de la forme
\begin{equation}\label{eqn:euler}
(\mathcal{L}-k\mathcal{I}) c^{(k)}=\mbox{polyn\^{o}me en}\quad
c_{i}^{( 0)},...,c_{i}^{(k-1)},\quad 1\leq i\leq n,\quad k\geq 1,
\end{equation}
o\`{u} $c^{(k)}=(c_1^{(k)},...,c_n^{(k)})^\top$ et
$$\mathcal{L}\equiv(\frac{\partial f_{i}}{\partial
w_{j}}(c_{1}^{(0)},...,c_{n}^{(0)}) +\delta _{ij}s_{i}) _{1\leq
i,j\leq n},$$ est la matrice jacobienne de (5). En outre, la
s\'{e}rie $(4) $ est convergente.
\end{thm}
\emph{Démonstration}: En substituant $\left(4\right) $ dans
$\left( 2\right) $, tout en tenant compte de la quasi-homogénéité
du système, on obtient
\begin{eqnarray}
\sum_{k=0}^{\infty }(k-s_{i}) c_{i}^{(k) }z^{k-s_{i}-1}&=&f_{i}(
\sum_{k=0}^{\infty }c_{1}^{(k) }z^{k-s_{1}},...,\sum_{k=0}^{\infty
}c_{n}^{( k)}z^{k-s_{n}})
,\nonumber\\
&=&f_{i}(z^{-s_{1}}(c_{1}^{(0) }+\sum_{k=1}^{\infty }c_{1}^{(
k)}z^{k}) ,...,z^{-s_{n}}(c_{n}^{(0)}+\sum_{k=1}^{\infty
}c_{n}^{(k)}z^{k})) ,\nonumber\\
&=&z^{-s_{i}-1}f_{i}(c_{1}^{(0)}+\sum_{k=1}^{\infty }c_{1}^{(k)
}z^{k},...,c_{n}^{(0)}+\sum_{k=1}^{\infty }c_{n}^{(k)
}z^{k}).\nonumber
\end{eqnarray}
Ensuite, on développe le second membre comme suit
\begin{eqnarray}
\sum_{k=0}^{\infty}(k-s_{i})
c_{i}^{(k)}z^{k}&=&f_{i}(c_{1}^{(0)},...,c_{n}^{(0)})+
\sum_{j=1}^{n}\frac{\partial f_{i}}{\partial
w_{j}}(c_{1}^{(0)},...,c_{n}^{(0)})
\sum_{k=1}^{\infty }c_{j}^{(k)}z^{k}\nonumber\\
&&+\sum_{k=2}^\infty z^k\sum_{(\alpha, \tau)\in
D_k}\frac{1}{\alpha !}\frac{\partial ^{\alpha }f_{i}}{\partial
w^{\alpha }}(c_{1}^{(0) },...,c_{n}^{( 0)})
\prod_{j=1}^{n}(c_{j}^{(\tau _{j})}) ^{\alpha _{j}},\nonumber
\end{eqnarray}
où $\alpha =\left( \alpha _{1},...,\alpha _{n}\right),$
$\tau=\left( \tau _{1},...,\tau _{n}\right) ,$
$$\left| \alpha \right| =\sum_{j=1}^{n}\alpha _{j},\quad
\alpha !=\prod_{j=1}^{n}\alpha _{j}!,$$
$$D_k=\{(\alpha, \tau):\tau_j>0, \forall j, \left|\alpha \right|>2, \sum_{j=1}^{n}\alpha
_{j}\tau _{j}=k\}.$$ En identifiant les termes ayant même
puissance au premier et au second membre, on obtient
successivement pour $k=0$ l'expression $\left(5\right),$ pour
$k=1,$ $$\left( \mathcal{L}-\mathcal{I}\right)c^{\left(1\right)
}=0,$$ et pour $k\geq 2,$
\begin{equation}\label{eqn:euler}
((\mathcal{L}-k\mathcal{I})c^{(k)})_i=-\sum_{(\alpha, \tau)\in
D_k}\frac{1}{\alpha !}\frac{\partial ^{\alpha }f_{i}}{\partial
w^{\alpha }}(c_{1}^{(0)},...,c_{n}^{(0)})
\prod_{j=1}^{n}(c_{j}^{(\tau _{j})}) ^{\alpha _{j}},
\end{equation}
où $\tau _{j}> 0,$ $\displaystyle{\sum_{j=1}^{n}\alpha _{j}\tau
_{j}=k},$ ce qui conduit aux expressions (explicites) $\left(
6\right) .$ La solution obtenue par la m\'{e}thode des
coefficients ind\'{e}termin\'{e}s est formelle du fait que nous
l'obtenons en effectuant sur des s\'{e}ries, que nous supposons a
priori convergentes, diverses op\'{e}rations dont la validit\'{e}
reste \`{a} justifier. Le th\'{e}or\`{e}me se trouvera donc
\'{e}tabli d\`{e}s que nous aurons v\'{e}rifi\'{e} que ces
s\'{e}ries sont convergentes. On utilise à cette fin la méthode
des fonctions majorantes ainsi que les travaux de M. Adler-P. van
Moerbeke [1] et J.P. Françoise [4]. Notons tout d'abord que des
param\`{e}tres libres apparaissent soit dans le syst\`{e}me
$\left( 5\right) $ de $n$ \'{e}quations \`{a} $n$ inconnues,
lorsque celui-ci admet un ensemble continue de solutions, soit par
le fait que $\lambda _{i}\equiv k\in \mathbb{N}^{*},1\leq i\leq
n,$ est une valeur propre de la matrice $\mathcal{L}.$ D\`{e}s
lors, les coefficients peuvent \^{e}tre vus comme \'{e}tant des
fonctions rationnelles sur une vari\'{e}t\'{e} affine $V,$ de
fibre le lieu
$$\bigcap_{i=1}^{n}\left\{ s_{i}c_{i}^{\left( 0\right) }+
f_{i}\left(c_{1}^{\left( 0\right) },...,c_{n}^{\left( 0\right)
}\right) =0\right\}.$$ Soit $n_0\in V$ et soit $K$ un
sous-ensemble compact de $V,$ contenant un voisinage ouvert de
$n_0.$ Notons que $K$ peut-être muni de la topologie du plan
complexe. Posons
$$A=1+\max \left\{\left| c_{1}^{\left(\tau _{1}\right)}(n_0)\right|,
\left| c_{2}^{\left( \tau _{2}\right) }(n_0)\right| ,...,\left|
c_{n}^{\left( \tau _{n}\right)}\right|(n_0) \right\}, 1\leq i\leq
n, 1\leq \tau _{i}\leq \lambda _{n},$$ où $\lambda_n$ désigne la
plus grande valeur propre de la matrice $\mathcal{L}.$ Soient $B$
et $C$ deux constantes avec $C>A$ telles que dans le compact $K$
on ait
$$\left| \frac{\partial ^{\alpha }f_{i}}{\partial w^{\alpha }}
(n_0) \right| \leq \alpha !B^{\left| \alpha \right| },$$$$\left|
\left( \mathcal{L}(n_0)-k\mathcal{I}_n\right) ^{-1}\right| \leq
C,\quad k\geq \lambda _{n}+1.$$ De $\left(7\right) $ on d\'{e}duit
que
$$\left| c_{i}^{\left( k\right)}(n_0)\right| \leq C\sum_{(\alpha, \tau)\in D}
B^{\left| \alpha \right| }\prod_{j=1}^{n}\left| c_{j}^{\left( \tau
_{j}\right) }\right| ^{\alpha _{j}},\quad k\geq \lambda _{n}+1.$$
Consid\'{e}rons maintenant la s\'{e}rie
$$\Phi\left( z\right) =Az+\sum_{k=2}^{\infty }\beta _{k}z^{k},$$
o\`{u} $\beta _{k}$ sont des nombres r\'{e}els d\'{e}finis
inductivement par $\beta_1\equiv A$ et
$$\beta _{k}\equiv C\sum_{(\alpha, \tau)\in D}B^{\left| \alpha
\right| }\prod_{j=1}^{n}\beta ^{\alpha _{j}}_{\tau _{j}},\quad
k\geq 2. $$ On v\'{e}rifie ais\'{e}ment par r\'{e}currence que la
s\'{e}rie $\Phi\left( z\right)$ est une majorante pour
$$\sum_{k=1}^{\infty }c_{i}^{\left( k\right) }z^{k},\quad 1\leq
i\leq n.$$ En effet, on a $\left| c_{i}^{\left( 1\right) }\right|
\leq A.$ Supposons que $\left| c_{i}^{\left( j\right) }\right|
\leq \beta _{j}, j < k, \forall i.$ Alors
\begin{eqnarray}
\left| c_{i}^{\left( k\right) }(n_0)\right|& \leq &
C\sum_{(\alpha, \tau)\in D}B^{\left| \alpha \right|
}\prod_{j=1}^{n}\left| c_{j}^{\left( \tau _{j}\right) }\right|
^{\alpha _{j}},\text{ }k\geq \lambda _{n}+1,\nonumber\\
&\leq &C\sum_{(\alpha, \tau)\in D}B^{\left| \alpha \right|
}\prod_{j=1}^{n}\left| \beta ^{\alpha _{j}}_{\tau _{j}}\right|,\nonumber\\
&=&\beta _{k}.\nonumber
\end{eqnarray}
D'autre part, il résulte de la d\'{e}finition des nombres $\beta
_{k}$ que $$\Phi \left( z\right) =Az+CB^{2}\frac{\left( n\Phi
\left( z\right) \right) ^{2}}{1-Bn\Phi \left( z\right) }.$$La
racine
$$
\Phi (z)=\frac{1+nABz-\sqrt{(1-2nAB(1+2nBC)
z+n^{2}A^{2}B^{2}z^{2})}}{2nB(1+nBC)},
$$ fournit la majorante cherch\'{e}e. D'o\`{u} la possibilit\'{e} d'un
d\'{e}veloppement en série entière au voisinage de $z=0.$
 Ceci achève la démonstration. $\square$ \\

\begin{rem}
La série (4) est l'unique solution méromorphe dans le sens où
cette solution r\'{e}sulte de ce que les
coefficients$\quad$$c_{i}^{\left( k\right) }$ se trouvent
d\'{e}termin\'{e}s de fa\c{c}on univoque avec la m\'{e}thode de
calcul adopt\'{e}.
\end{rem}
\begin{rem}
Le résultat du théorème précédent s'applique à l'équation
différentielle quasi-homogène d'ordre $n$ suivante :
\begin{eqnarray}\label{eqn:euler}
\frac{d^{n}w}{dz^{n}}=f\left(w,\frac{dw}{dz},...,\frac{d^{n-1}w}{dz^{n-1}}\right).
\end{eqnarray}
$f$ étant une fonction rationnelle en
$w,\frac{dw}{dz},...,\frac{d^{n-1}w}{dz^{n-1}}$ et
\begin{eqnarray}
w(z_{0})&=&w_{1}^{0},\nonumber\\
\frac{dw}{dz}(z_{0})&=&w_{2}^{0},\nonumber\\
&\vdots&\nonumber\\
\frac{d^{n-1}w}{dz^{n-1}}(z_{0})&=&w_{n}^{0}.\nonumber
\end{eqnarray}
En effet, l'équation $(8)$ se ramène à un système de $n$ équations
du premier ordre en posant
\begin{eqnarray}
w\left( z\right)&=&w_{1}\left( z\right),\nonumber\\
\frac{dw}{dz}\left( z\right)&=&w_{2}\left(
z\right),\nonumber\\
&\vdots&\nonumber\\
\frac{d^{n-1}w}{dz^{n-1}}\left( z\right)&=&w_{n}\left(
z\right).\nonumber
\end{eqnarray}
On obtient ainsi
\begin{eqnarray}
\frac{dw_{1}}{dz}&=&w_{2},\nonumber\\
\frac{dw_{2}}{dz}&=&w_{3},\nonumber\\
&\vdots&\nonumber\\
\frac{dw_{n-1}}{dz}&=&w_{n},\nonumber\\
\frac{dw_{n}}{dz}&=&f\left(w_{1},w_{2},...,w_{n}\right).\nonumber
\end{eqnarray}
Un tel syst\`{e}me constitue un cas particulier du syst\`{e}me
normal $\left(2\right)$.
\end{rem}

\section{Connection avec l'int\'{e}grabilit\'{e} alg\'{e}brique des systèmes hamiltoniens}

Considérons un sytème hamiltonien
\begin{equation}\label{eqn:euler}
\frac{dw}{dz}=J\left(w\right) \frac{\partial H}{\partial w},\quad
w\in \mathbb{R}^{n},\quad n=2m+k, \end{equation} o\`{u} $H$ est
l'hamiltonien et $J\left( w\right) $ est une matrice r\'{e}elle
antisym\'{e}trique telle que les crochets de Poisson
correspondants v\'{e}rifient l'identit\'{e} de Jacobi:$$\left\{
\left\{ H,F\right\} ,G\right\} +\left\{ \left\{ F,G\right\}
,H\right\} +\left\{ \left\{ G,H\right\} ,F\right\} =0, \forall H,
F, G\in \mathcal{C}^\infty (\mathbb{R}^n),$$ o\`{u}
$$
\left\{ H,F\right\}=\sum_{i,j}J_{ij}\frac{\partial H}{\partial
w_{i}}\frac{\partial F}{\partial w_{j}}.$$ On suppose que le
syst\`{e}me $(9)$ est complètement intégrable c'est-à-dire qu'il
admet $m+k$ intégrales premières $H_{1}=H,H_{2},...,H_{m+k}$
fonctionnellement indépendantes dont $m$ intégrales sont en
involution (i.e.,\quad $\{ H_{i},H_{j}\} =0,$ $1\leq i,j\leq m,$),
$k$ int\'{e}grales sont des fonctions de Casimir (i.e.,\quad
$J\frac{\partial H_{m+i}}{\partial w}=0,$ $1\leq i\leq k$) et
telles que pour presque tous les $c_{i}\in \mathbb{R}$ les
vari\'{e}t\'{e}s invariantes
\begin{equation}\label{eqn:euler}
\overset{m+k}{\underset{i=1}{\bigcap}}\{w\in
\mathbb{R}^{n}:H_{i}(w) =c_{i}\} ,
\end{equation}
sont compactes et connexes. D'apr\`{e}s le th\'{e}or\`{e}me
d'Arnold-Liouville $\left[ 10\right]$, les vari\'{e}t\'{e}s
$\left( 10\right) $ sont diff\'{e}omorphes aux tores r\'{e}els
$T_{\mathbb{R}}^{m}=\mathbb{R}^{m}/r\acute{e}seau.$ En outre les
flots d\'{e}finis par les champs de vecteurs $W_{H_{i}},1\leq
i\leq m,$ sont des mouvements rectilignes sur ce tore et les
\'{e}quations du probl\`{e}me sont int\'{e}grables par
quadratures.\\
Soient maintenant $w\in \mathbb{C}^{n}$, $z\in \mathbb{C}$ et
$\Delta \subset $ $\mathbb{C}^{n}$ un ouvert non vide de Zariski.
Comme $H_{1},...,H_{m+k}$ sont fonctionnellement
ind\'{e}pendantes, alors l'application
$$\varphi =\left( H_{1},...,H_{m+k}\right) :\mathbb{C}^{n}\longrightarrow
\mathbb{C}^{m+k},$$est une submersion g\'{e}n\'{e}rique sur
$\Delta .$ Soit $ \textbf{I}=\varphi (\mathbb{C}^{n}\backslash
\Delta),$ le lieu critique de $\varphi $ et d\'{e}signons par $adh
\textbf{I}$ l'adh\'{e}rence (ou fermeture) de Zariski de
$\textbf{I}$ dans $\mathbb{C}^{m+k}.$

\begin{defn}
Le syst\`{e}me différentiel $\left(9\right)$ dont le c\^{o}t\'{e}
droit est polynomial est alg\'{e}briquement compl\`{e}tement
int\'{e}grable si pour $c=(c_1,...,c_{m+k})\in \mathbb{C}^{m+k}$
$\backslash $ $adh \textbf{I},$ la fibre
\begin{equation}\label{eqn:euler}
M_c\equiv\varphi^{-1}\left( c\right) =\bigcap_{i=1}^{m+k}\left\{
w\in \mathbb{C}^{n}:H_{i}\left( w\right) =c_{i}\right\} ,
\end{equation}
est la partie affine d'une vari\'{e}t\'{e} ab\'{e}lienne (i.e., un
tore complexe $T_{\mathbb{C}}^{m}\simeq
\mathbb{C}^{m}/r\acute{e}seau$ qui poss\`{e}de un plongement dans
un espace projectif ). En outre, les flots $g_{W_{i}}^{z}\left(
w\right) ,$ $w\in M_c, $ $z\in \mathbb{C},$ d\'{e}finies par les
champs de vecteurs $W_{H_{1}},...,W_{H_{m}}$ sont des lignes
droites sur $T_{\mathbb{C}}^{m}$ c'est-\`{a}-dire $$\left[
g_{W_{i}}^{z}\left( w\right) \right] _{j}=f_{j}\left(
p+z(k_{1}^{i},...,k_{n}^{i})\right) ,$$o\`{u} $f_{j}\left(
z_{1},...,z_{m}\right) $ sont des fonctions ab\'{e}liennes
(m\'{e}romorphes) sur le tore $T_{\mathbb{C}}^{m},$
$f_{j}(p)=w_{j},$ $1\leq j\leq n.$
\end{defn}

Soit $\overline{M_c}$\ la fermeture projective\ de $M_c$ dans
l'espace projectif complexe $\mathbb{CP}^{n}$\ de dimension $n.$
Alors $\overline{M_c}$ n'est pas une vari\'{e}t\'{e} ab\'{e}lienne
puisque cette derni\`{e}re n'est pas simplement connexe et ne peut
donc en général \^{e}tre une intersection compl\`{e}te projective.
D\`{e}s lors, pour que $M_c$ soit la partie affine d'une
vari\'{e}t\'{e} ab\'{e}lienne, la vari\'{e}t\'{e} $\overline{M_c}$
doit \^{e}tre singuli\`{e}re \`{a} l'infini. En \'{e}clatant la
singularit\'{e} le long du lieu atteint par le flot et en
implosant la partie du lieu qui n'est pas atteint par le flot, on
montre que la vari\'{e}t\'{e} $\overline{M_c}$ se transforme en
une vari\'{e}t\'{e} ab\'{e}lienne $\widetilde{M_c}$ et le lieu
\`{a} l'infini se transforme en une ou plusieurs
sous-vari\'{e}t\'{e}s de codimension $1.$ Et c'est l\`{a} o\`{u}
le th\'{e}or\`{e}me 1, va jouer un r\^{o}le crucial. On
proc\`{e}de comme suit : soit $w_{i}\longrightarrow u_{i}$ une
transformation birationnelle telle qu'au voisinage de $z=0,$ on
ait :
\begin{eqnarray}
u_{i}&=&\alpha _{i}+\circ\left( z\right) ,\text{ }1\leq i\leq
n-1,\nonumber\\
u_{n}&=&z+\circ\left( z^{2}\right) ,\nonumber
\end{eqnarray}
o\`{u} $\alpha _{1},...,\alpha _{n-1}$ sont des param\`{e}tres
libres. Les nouvelles variables $u_{i}$ ont pour effet
d'\'{e}clater la singularit\'{e} de la vari\'{e}t\'{e} projective
$\overline{M_c}$ le long du lieu \`{a} l'infini atteint par le
flot. Exprim\'{e}es dans ces nouvelles variables
$u_{1},...,u_{n},$ les \'{e}quations diff\'{e}rentielles sont
r\'{e}guli\`{e}res et holomorphes au voisinage de $u_{n}=0$ tandis
que les \'{e}quations d\'{e}finissant la fibre $M_c$
s'\'{e}crivent sous la forme :
$$F_{i}\left( u_{1}\left( z\right) ,...,u_{n-1}\left( z\right)
,u_{n}\left( z\right) \right) =c_{i},\text{ }1\leq i\leq m+k,$$
o\`{u} $F_{1},...,F_{m+k}$ sont des polyn\^{o}mes en $w$. Pour
$z=0,$ on obtient $$F_{i}\left( \alpha _{1},...,\alpha
_{n-1},0\right) =a_{i},\text{ }1\leq i\leq m+k,$$ et ces relations
alg\'{e}briques entre les param\`{e}tres libres $\alpha
_{1},...,\alpha _{n-1}$ fournissent les \'{e}quations d'une
sous-vari\'{e}t\'{e} $\mathcal{D}$ qui jouera, entre autres, un
rôle important dans la compactification de la fibre $M_c.$ En
fait, les param\`{e}tres libres $\alpha _{1},...,\alpha _{n-1}$ et
la sous vari\'{e}t\'{e} $\mathcal{D}$ peuvent s'obtenir
directement de la mani\`{e}re suivante: d'abord l'on montre
l'existence de solutions $w=(w_{1},w_{2},\ldots ,w_{n})$ du
syst\`{e}me $\left( 9\right) $ sous la forme de s\'{e}ries de
Laurent $\left( 4\right) $ d\'{e}pendant de $n-1$ param\`{e}tres
libres $\alpha _{1},...,\alpha _{n-1}.$ En substituant ces
d\'{e}veloppements dans le syst\`{e}me $\left( 9\right) $, on voit
(d'apr\`{e}s le th\'{e}or\`{e}me 1) que les coefficients
$c^{\left( 0\right) },c^{\left( 1\right) },...,$ satisfont aux
\'{e}quations $\left( 5\right) $ et $\left( 6\right).$ L'\'{e}tape
suivante consiste \`{a} consid\'{e}rer la fermeture $\mathcal{D}$
des composantes continues de l'ensemble des s\'{e}ries de Laurent
de $w\left( z\right) $ tels que: $H_{1}\left( w\right)
=a_{1},\ldots ,$ $H_{m+k}\left( w\right) =a_{m+k}.$ Plus
pr\'{e}cisement,$$\mathcal{D}=\bigcap_{i=1}^{m+k}\left\{
\text{coefficient de }z^{0}\text{ dans }H_{i}\left( w\left(
z\right) \right) =a_{i}\right\} .$$ C'est une sous-vari\'{e}t\'{e}
(un diviseur) de codimension $1.$ Ensuite on proc\`{e}de \`{a} la
compactification de la fibre $M_c$ $\left( 11\right) $ en une
vari\'{e}t\'{e} ab\'{e}lienne $\widetilde{M_c}$. Cette
compactification s'obtient par l'adjonction \`{a} $M_c$ de ce
diviseur $\mathcal{D}$.

\subsection{Rotation d'un corps solide autour d'un point fixe}

Les \'{e}quations du mouvement d'un corps solide autour d'un point
fixe s'\'{e}crivent sous la forme
\begin{eqnarray}
\frac{dM}{dt}&=&M\wedge \Omega +\mu g\text{ }\Gamma \wedge L,\\
\frac{d\Gamma }{dt}&=&\Gamma \wedge \Omega ,\nonumber
\end{eqnarray}
o\`{u} $\wedge $ est le produit vectoriel dans $\mathbb{R}^{3},$
$M=\left( m_{1},m_{2},m_{3}\right) $ le moment angulaire du
solide$,\Omega =\left( m_{1}/I_{1},m_{2}/I_{2},m_{3}/I_{3}\right)
$ la vitesse angulaire, $I_{1},I_{2}$ et $I_{3},$ les moments
d'inertie, $\Gamma =\left( \gamma _{1},\gamma _{2},\gamma
_{3}\right) $ le vecteur vertical unitaire, $\mu $ la masse du
solide, $g$ l'acc\'{e}l\'{e}ration de la pesanteur, et enfin,
$L=\left( l_{1},l_{2},l_{3}\right) $ le vecteur unitaire ayant
pour origine le point fixe et dirig\'{e} vers le centre de
gravit\'{e}; tous ces vecteurs sont consid\'{e}r\'{e}s dans un
syst\`{e}me mobile dont les coordonn\'{e}es sont fix\'{e}es aux
axes principaux d'inertie. L'espace de configuration d'un solide
avec un point fixe est le groupe des rotations:
$$SO\left( 3\right) =\left\{ U\text{ matrice d'ordre trois}:\text{ }U^{-1}=
U^{\top },\text{ }\det U=1\right\}.$$ C'est le groupe des matrices
orthogonales d'ordre trois et le mouvement de ce solide est
d\'{e}crit par une courbe sur ce groupe. L'espace des vitesses
angulaires de toutes les rotations est l'alg\`{e}bre de Lie du
groupe $SO\left( 3\right) ;$ c'est l'alg\`{e}bre
$$so(3)=\left\{ A\text{ matrice d'ordre trois}:U^{\top }+U=0\right\}
,$$ des matrices antisym\'{e}triques d'ordre trois$.$ Cette
alg\`{e}bre est engendr\'{e}e comme espace vectoriel par les
matrices
$$
{e_{1}}=\left(\begin{array}{ccc}
0&0&0\\
0&0&-1\\
0&1&0
\end{array}\right),\qquad
{e_{2}}=\left(\begin{array}{ccc}
0&0&1\\
0&0&0\\
-1&0&0
\end{array}\right),\qquad
{e_{1}}=\left(\begin{array}{ccc}
0&-1&0\\
1&0&0\\
0&0&0
\end{array}\right),
$$
qui v\'{e}rifient les relations de commutation $$\left[
e_{1},e_{2}\right] =e_{3},\quad\left[ e_{2},e_{3}\right]
=e_{1},\quad\left[ e_{3},e_{1}\right] =e_{2}.$$ On utilisera dans
la suite le fait que si l'on identifie $so\left( 3\right) $ \`{a}
$\mathbb{R}^{3}$ en envoyant $\left( e_{1},e_{2},e_{3}\right) $
sur la base canonique de $\mathbb{R}^{3},$ le crochet de $so\left(
3\right) $ correspond au produit vectoriel. En d'autres termes,
consid\'{e}rons l'application
$$\mathbb{R}^{3}\longrightarrow so(3),\text{ }a=
\left( a_{1},a_{2},a_{3}\right) \longmapsto
{A}=\left(\begin{array}{ccc}
0&-a_{3}&a_{2}\\
a_{3}&0&-a_{1}\\
-a_{2}&a_{1}&0
\end{array}\right),
$$
laquelle d\'{e}finit un isomorphisme entre les alg\`{e}bres de Lie
$\left( \mathbb{R}^{3},\wedge \right) $ et $\left( so(3),\left[
,\right] \right) $ o\`{u}
$$a\wedge b\longmapsto \left[ A,B\right] =AB-BA.$$
En utilisant cet isomorphisme, on peut r\'{e}ecrire le système
$\left(12\right) $ sous la forme
\begin{eqnarray}
\frac{dM}{dt}&=&\left[ M,\Omega \right] +\mu g\text{ }\left[
\Gamma,L\right] ,\nonumber\\
\frac{d\Gamma }{dt}&=&\left[ \Gamma ,\Omega \right] ,\nonumber
\end{eqnarray}
o\`{u}
$$M=\left( M_{ij}\right) _{1\leq i,j\leq 3}\equiv
\sum_{i=1}^{3}m_{i}e_{i}\equiv \left(\begin{array}{ccc}
0&-m_{3}&m_{2}\\
m_{3}&0&-m_{1}\\
-m_{2}&m_{1}&0
\end{array}\right)\in so\left( 3\right) ,
$$
$$\Omega =\left( \Omega _{ij}\right) _{1\leq i,j\leq 3}\equiv
\sum_{i=1}^{3}\omega _{i}e_{i}\equiv \left(\begin{array}{ccc}
0&-\omega_{3}&\omega_{2}\\
\omega_{3}&0&-\omega_{1}\\
-\omega_{2}&\omega_{1}&0
\end{array}\right)\in so\left( 3\right) ,
$$
$$\Gamma =\left( \gamma _{ij}\right) _{1\leq i,j\leq 3}\equiv
\sum_{i=1}^{3}\gamma _{i}e_{i}\equiv \left(\begin{array}{ccc}
0&-\gamma_{3}&\gamma_{2}\\
\gamma_{3}&0&-\gamma_{1}\\
-\gamma_{2}&\gamma_{1}&0
\end{array}\right)\in so\left( 3\right) ,
$$
et
$$L = \left(\begin{array}{ccc}
0&-l_{3}&l_{2}\\
l_{3}&0&-l_{1}\\
-l_{2}&l_{1}&0
\end{array}\right)\in so\left( 3\right) ,
$$
En tenant compte du fait que $M=I\Omega ,$ alors les \'{e}quations
p\'{e}c\'{e}dentes deviennent
\begin{eqnarray}
\frac{dM}{dt}&=&\left[ M,\Lambda M\right] +\mu g\text{ }\left[
\Gamma ,L\right] ,\\
\frac{d\Gamma }{dt}&=&\left[ \Gamma ,\Lambda M\right] ,\nonumber
\end{eqnarray}
où
$$\Lambda M=\left( \Lambda _{ij}M_{ij}\right) _{1\leq i,j\leq 3}\equiv
\sum_{i=1}^{3}\lambda _{i}m_{i}e_{i}\equiv
\left(\begin{array}{ccc}
0&-\lambda _{3}m_{3}&\lambda _{2}m_{2}\\
\lambda _{3}m_{3}&0&-\lambda _{1}m_{1}\\
-\lambda _{2}m_{2}&\lambda _{1}m_{1}&0
\end{array}\right)\in so\left( 3\right),
$$
avec $\lambda _{i}\equiv I_{i}^{-1}.$ Le syst\`{e}me $\left(
13\right) $ est compl\`{e}tement int\'{e}grable seulement dans les
cas suivants :

$a)$ \underline{Cas d'Euler} : Dans ce cas, on a
$$l_{1}=l_{2}=l_{3}=0,$$
c'est-\`{a}-dire le point fixe est son centre de gravit\'{e}.
Autrement dit, Les équations d'Euler (On parle aussi de mouvement
d'Euler-Poinsot du solide) du mouvement de rotation d'un solide
autour d'un point fixe, pris comme origine du repère lié au
solide, lorsqu'aucune force extérieure n'est appliquée au système,
peuvent s'écrire forme explicite
\begin{eqnarray}
\frac{dm_{1}}{dt}&=&\left( \lambda _{3}-\lambda _{2}\right)
m_{2}m_{3},\nonumber\\
\frac{dm_{2}}{dt}&=&\left( \lambda _{1}-\lambda _{3}\right) m_{1}m_{3},\\
\frac{dm_{3}}{dt}&=&\left( \lambda _{2}-\lambda _{1}\right)
m_{1}m_{2}.\nonumber
\end{eqnarray}
Ces équations forment un champ de vecteurs hamiltonien de la forme
(9) avec $n=3$, $m=k=1$, $z=t$, $w= \left(
m_{1},m_{2},m_{3}\right)^{\intercal}$,
$$H=\frac{1}{2}\left( \lambda _{1}m_{1}^{2}+\lambda
_{2}m_{2}^{2}+\lambda _{3}m_{3}^{2}\right),$$ l'hamiltonien et
$$J=\left(\begin{array}{ccc}
0&-m_{3}&m_{2}\\
m_{3}&0&-m_{1}\\
-m_{2}&m_{1}&0
\end{array}\right)\in so\left( 3\right).
$$
Ces équations admettent deux intégrales premières quadratiques :
$H_1=H$ et
$$H_2=\frac{1}{2}\left( m_{1}^{2}+m_{2}^{2}+m_{3}^{2}\right).$$
Ces intégrales sont fonctionnellement indépendantes, en involution
et le système en question est complètement intégrable. La
r\'{e}solution explicite des équations  d'Euler est délicate dans
le cas général où $\lambda _{1}$, $\lambda _{2}$ et $\lambda _{3}$
sont tous différents; les solutions s'expriment à l'aide de
fonctions elliptiques de Jacobi comme suit (pour le détail voir
par exemple [15]) :
\begin{equation}\label{eqn:euler}
\left\{\begin{array}{rl}
m_1=\sqrt{\frac{2H_1-H_2\lambda_{3}}{\lambda_1-\lambda_3}}&\mathbf{cn}(t\sqrt{(\lambda
_{2}-\lambda _{3})(H_2\lambda
_{1}-2H_1)}),\\
m_2=\sqrt{\frac{2H_1-H_2\lambda_{3}}{\lambda_2-\lambda_3}}&\mathbf{sn}(t\sqrt{(\lambda
_{2}-\lambda _{3})(H_2\lambda
_{1}-2H_1)}),\\
m_3=\sqrt{\frac{H_2\lambda_{1}-2H_1}{\lambda_1-\lambda_3}}&\mathbf{dn}(t\sqrt{(\lambda
_{2}-\lambda _{3})(H_2\lambda _{1}-2H_1)}).
\end{array}\right.
\end{equation}
Le mouvement de $\left( m_{1},m_{2},m_{3}\right) $ s'effectue sur
l'intersection $M_c$ d'un ellipsoide avec une sphère. Les deux
cercles d\'{e}finies par $M_c,$ forme la partie r\'{e}elle d'un
tore complexe de dimension $1,$ d\'{e}finie par la courbe
elliptique $\mathcal{E}$ :
$$
\mathcal{E}: y^2=(1- s^{2})(1-k^2 s^{2}).
$$
où
$$k^2\equiv\frac{(\lambda _{1}-\lambda
_{2})(2H_1-H_2\lambda _{3})}{(\lambda_2-\lambda_3)(H_2\lambda
_{1}-2H_1)}.$$ et
$$
s\equiv m_2\sqrt{\frac{\lambda _{2}-\lambda _{3}}{2H_1-H_2\lambda
_{3}}}
$$
L'intersection complexe $\left( \subset \mathbb{C}^{3}\right) $
est la partie affine d'une courbe elliptique $\overline{M_c}\in
\mathbb{CP}^{3}$. On montre que $\overline{M_c}$ est isomorphe
\`{a} la courbe elliptique $\mathcal{E}$. En outre l'intersection
réelle $\left( \subset \mathbb{R}^{3}\right) $ s'\'{e}tend au tore
complexe $\mathbb{C}/r\acute{e}seau$ et le flot se lin\'{e}arise
sur ce tore. Si $p\left( t\right) =\left( m_{1}\left( t\right)
,m_{2}\left( t\right) ,m_{3}\left( t\right) \right)$, est une
solution de $\left( 14\right)$, la loi reliant $p\left(
t_{1}+t_{2}\right)$ \`{a} $p\left( t_{1}\right) $ et $p\left(
t_{2}\right)$ est la loi d'addition sur la courbe elliptique (voir
[15]). D'apr\`{e}s les \'{e}quations $\left( 14\right) $, l'unique
diff\'{e}rentielle holomorphe sur $\overline{M_c}$ est donn\'{e}e
par $$\omega =\frac{dm_{1}}{\left( \lambda _{3}-\lambda
_{2}\right) m_{2}m_{3}}=\frac{dm_{2}}{\left( \lambda _{1}-\lambda
_{3}\right) m_{1}m_{3}}=\frac{dm_{3}}{\left( \lambda _{2}-\lambda
_{1}\right) m_{1}m_{2}},$$d'o\`{u}$$t=\int_{p\left( 0\right)
}^{p\left( t\right) }\omega ,\text{ \quad }p\left( 0\right) \in
\overline{M_c}.$$ Le syst\`{e}me $\left( 14\right) $ est invariant
par les transformations $$t\longrightarrow\alpha ^{-1}t,\text{
}m_{1}\longrightarrow \alpha m_{1},\text{ }m_{2}\longrightarrow
\alpha m_{2},\text{ }m_{3}\longrightarrow \alpha m_{3}.$$
Celles-ci sont uniques puisque le déterminant $\left( 3\right) $
est \'{e}gal \`{a}
\begin{eqnarray}
\Delta&=&\left|
\begin{array}{ccc}
-\left( \lambda _{3}-\lambda _{2}\right) m_{2}m_{3}&\left( \lambda
_{3}-\lambda _{2}\right) m_{2}m_{3}&\left( \lambda _{3}-\lambda
_{2}\right) m_{2}m_{3}\\
\left( \lambda _{1}-\lambda _{3}\right) m_{1}m_{3}&-\left( \lambda
_{1}-\lambda _{3}\right) m_{1}m_{3}&\left( \lambda _{1}-\lambda
_{3}\right) m_{1}m_{3}\\
\left( \lambda _{2}-\lambda _{1}\right) m_{1}m_{2}&\left( \lambda
_{2}-\lambda _{1}\right) m_{1}m_{2}&-\left( \lambda _{2}-\lambda
_{1}\right) m_{1}m_{2}
\end{array}
\right|,\nonumber\\
&=&4\left( \lambda _{3}-\lambda _{2}\right) \left( \lambda
_{1}-\lambda _{3}\right) \left( \lambda _{2}-\lambda _{1}\right)
m_{1}^{2}m_{2}^{2}m_{3}^{2},\nonumber\\
&\neq &0.\nonumber
\end{eqnarray}
On peut donc chercher des solutions du syst\`{e}me $\left(
14\right)$ sous la forme de s\'{e}ries de Laurent
\begin{eqnarray}
m_{1}&=&\frac{1}{t}\left( a_{0}+a_{1}t+a_{2}t^{2}+\cdots
\right),\nonumber\\
m_{2}&=&\frac{1}{t}\left( b_{0}+b_{1}t+b_{2}t^{2}+\cdots \right),\nonumber\\
m_{3}&=&\frac{1}{t}\left( c_{0}+c_{1}t+c_{2}t^{2}+\cdots
\right),\nonumber
\end{eqnarray}
d\'{e}pendant de $\dim (\mbox{espace de phase})-1=2$
param\`{e}tres libres. En substituant ces équations dans le
système $\left( 14\right) $, on voit que :

1) les coefficients $a_0,b_0,c_0,$ satisfont aux \'{e}quations
\begin{eqnarray}
a_{0}+\left( \lambda _{3}-\lambda _{2}\right) b_{0}c_{0}&=&0,\nonumber\\
b_{0}+\left( \lambda _{1}-\lambda _{3}\right) a_{0}c_{0}&=&0,\nonumber\\
c_{0}+\left(\lambda_{2}-\lambda_{1}\right)a_{0}b_{0}&=&0,\nonumber
\end{eqnarray}
dont les solutions sont\\
\underline{$1^{er}$ cas} : $a_{0}=\frac{-1}{\sqrt{\left( \lambda
_{2}-\lambda _{1}\right) \left( \lambda _{1}-\lambda _{3}\right)
}},\quad b_{0}=\frac{1}{\sqrt{\left( \lambda _{2}-\lambda
_{1}\right) \left( \lambda _{3}-\lambda _{2}\right) }},\quad
c_{0}=\frac{1}{\sqrt{\left( \lambda _{1}-\lambda _{3}\right)
\left( \lambda _{3}-\lambda _{2}\right) }}.$\\
\underline{$2^{ème}$ cas} : $a_{0}=\frac{1}{\sqrt{\left( \lambda
_{2}-\lambda _{1}\right) \left( \lambda _{1}-\lambda _{3}\right)
}},\quad b_{0}=\frac{1}{\sqrt{\left( \lambda _{2}-\lambda
_{1}\right) \left( \lambda _{3}-\lambda _{2}\right) }},\quad
c_{0}=\frac{-1}{\sqrt{\left( \lambda _{1}-\lambda
_{3}\right) \left( \lambda _{3}-\lambda _{2}\right) }}.$\\
\underline{$3^{ème}$ cas} : $a_{0}=\frac{1}{\sqrt{\left( \lambda
_{2}-\lambda _{1}\right) \left( \lambda _{1}-\lambda _{3}\right)
}},\quad b_{0}=\frac{-1}{\sqrt{\left( \lambda _{2}-\lambda
_{1}\right) \left( \lambda _{3}-\lambda _{2}\right) }},\quad
c_{0}=\frac{1}{\sqrt{\left( \lambda _{1}-\lambda _{3}\right)
\left( \lambda _{3}-\lambda _{2}\right) }}.$\\
\underline{$4^{ème}$ cas} : $a_{0}=\frac{-1}{\sqrt{\left( \lambda
_{2}-\lambda _{1}\right) \left( \lambda _{1}-\lambda _{3}\right)
}},\quad b_{0}=\frac{-1}{\sqrt{\left( \lambda _{2}-\lambda
_{1}\right) \left( \lambda _{3}-\lambda _{2}\right) }},\quad
c_{0}=\frac{-1}{\sqrt{\left( \lambda _{1}-\lambda _{3}\right)
\left( \lambda _{3}-\lambda _{2}\right) }}$.

2) les coefficients $a_1,b_1,c_1,$ satisfont aux \'{e}quations
\begin{eqnarray}
\left( \lambda _{3}-\allowbreak \lambda _{2}\right)
b_{0}c_{1}+\left( \lambda _{3}
-\lambda _{2}\right) b_{1}c_{0}&=&0,\nonumber\\
\left( \lambda _{1}-\allowbreak \lambda _{3}\right)
a_{0}c_{1}+\left( \lambda _{1}
-\lambda _{3}\right) a_{1}c_{0}&=&0,\nonumber\\
\left( \lambda _{2}-\allowbreak \lambda _{1}\right)
a_{0}b_{1}+\left( \lambda _{2}-\lambda _{1}\right)
a_{1}b_{0}&=&0,\nonumber
\end{eqnarray}
dont les solutions sont dant tous les cas : $a_{1}=b_{1}=c_{1}=0$.

3) les coefficients $a_2,b_2,c_2,$ satisfont aux \'{e}quations
\begin{eqnarray}
&&a_{2}-\lambda _{3}b_{0}c_{2}-\lambda _{3}b_{1}c_{1}-\lambda
_{3}b_{2}c_{0} +\allowbreak \lambda _{2}b_{0}c_{2}+\lambda
_{2}b_{1}c_{1}
+\allowbreak \lambda _{2}b_{2}c_{0}=0,\nonumber\\
&&b_{2}-\lambda _{1}a_{0}c_{2}-\lambda _{1}a_{1}c_{1}-\lambda
_{1}a_{2}c_{0} +\allowbreak \lambda _{3}a_{0}c_{2}+\lambda
_{3}a_{1}c_{1}
+\allowbreak \lambda _{3}a_{2}c_{0}=0,\nonumber\\
&&c_{2}-\lambda _{2}a_{0}b_{2}-\lambda _{2}a_{1}b_{1}-\lambda
_{2}a_{2}b_{0}+\allowbreak \lambda _{1}a_{0}b_{2}+\lambda
_{1}a_{1}b_{1}+\allowbreak \lambda _{1}a_{2}b_{0}=0,\nonumber
\end{eqnarray}
dont les solutions qui correspondent aux différents cas sont
respectivement :\\
\underline{$1^{er}$ cas} : $a_{2}=\frac{\sqrt{\left( \lambda
_{3}-\lambda _{2}\right) }}{\sqrt{\left( \lambda _{1}-\lambda
_{3}\right) }}b_{2}+\frac{\sqrt{\left( \lambda _{3}-\lambda
_{2}\right) }}{\sqrt{\left( \lambda _{2}-\lambda _{1}\right)
}}c_{2}$.\\
\underline{$2^{ème}$ cas} : $a_{2}=-\frac{\sqrt{\lambda
_{3}-\lambda _{2}}}{\sqrt{\lambda _{1}-\lambda
_{3}}}b_{2}+\frac{\sqrt{\lambda _{3}-\lambda _{2}}}{\sqrt{\lambda
_{2}-\lambda _{1}}}c_{2}$.\\
\underline{$3^{ème}$ cas} : $a_{2}=\frac{\sqrt{\lambda
_{3}-\lambda _{2}}}{\sqrt{\lambda _{1}-\lambda
_{3}}}b_{2}-\frac{\sqrt{\lambda _{3}-\lambda _{2}}}{\sqrt{\lambda
_{2}-\lambda _{1}}}c_{2}$.\\
\underline{$4^{ème}$ cas} : $a_{2}=-\frac{\sqrt{\lambda
_{3}-\lambda _{2}}}{\sqrt{\lambda _{1}-\lambda
_{3}}}b_{2}-\frac{\sqrt{\lambda _{3}-\lambda _{2}}}{\sqrt{\lambda
_{2}-\lambda _{1}}}c_{2}$.\\
où $b_2$ et $c_2$ sont deux paramètres libres.\\
Par conséquent, pour le premier cas on a
\begin{eqnarray}
m_{1}&=&\frac{-1}{t\sqrt{\left( \lambda _{2}-\lambda _{1}\right)
\left( \lambda _{1}-\lambda _{3}\right) }}+\left(
\frac{\sqrt{\left( \lambda _{3}-\lambda _{2}\right)
}}{\sqrt{\left( \lambda _{1}-\lambda _{3}\right)
}}b_{2}+\frac{\sqrt{\left( \lambda _{3}-\lambda _{2}\right)
}}{\sqrt{\left( \lambda _{2}-\lambda _{1}\right) }}c_{2}\right)
t+\cdots ,\nonumber\\
m_{2}&=&\frac{1}{t\sqrt{\left( \lambda _{2}-\lambda _{1}\right)
\left( \lambda _{3}-\lambda _{2}\right)
}}+b_{2}t+\cdots ,\nonumber\\
m_{3}&=&\frac{1}{t\sqrt{\left( \lambda _{1}-\lambda _{3}\right)
\left( \lambda _{3}-\lambda _{2}\right) }}+c_{2}t+\cdots.\nonumber
\end{eqnarray}
En substituant ces développements dans les intégrales premières
$H_1$ et $H_2$, on obtient
\begin{eqnarray}
H_{1}&=&2\frac{\sqrt{\lambda _{3}-\lambda _{2}}}{\sqrt{\lambda
_{2}-\lambda _{1}}}\left( \frac{1}{\lambda _{3}-\lambda
_{2}}-\frac{1}{\lambda _{1}-\lambda _{3}}\right)
b_{2}+2\frac{\sqrt{\lambda _{3}-\lambda _{2}}}{\sqrt{\lambda
_{1}-\lambda _{3}}}\left( \frac{1}{\lambda _{3}-\lambda
_{2}}-\allowbreak \frac{1}{\lambda _{2}-\lambda _{1}}\right)
c_{2},\nonumber\\
H_{2}&=&2\frac{\sqrt{\lambda _{3}-\lambda _{2}}}{\sqrt{\lambda
_{2}-\lambda _{1}}}\left( \frac{\lambda _{2}}{\lambda _{3}-\lambda
_{2}}-\frac{\lambda _{1}}{\lambda _{1}-\lambda _{3}}\right) b_{2}
+2\frac{\sqrt{\lambda _{3}-\lambda _{2}}}{\sqrt{\lambda
_{1}-\lambda _{3}}}\left( \frac{\lambda _{3}}{\lambda _{3}-\lambda
_{2}}-\allowbreak \frac{\lambda _{1}}{\lambda _{2}-\lambda
_{1}}\right) c_{2},\nonumber
\end{eqnarray}
et on en déduit les relations
\begin{eqnarray}
c_{2}&=&\frac{1}{6\sqrt{\left( \lambda _{1}-\lambda _{3}\right)
\left( \lambda _{3}-\lambda _{2}\right) }}\left( \left( \lambda
_{3}-\lambda _{2}\right) \left( \lambda _{1}H_{1}-\allowbreak
H_{2}\right) -\left( \lambda _{1}-\lambda _{3}\right) \left(
\lambda _{2}H_{1}-H_{2}\right) \right),\nonumber\\
b_{2}&=&\frac{1}{6\sqrt{\left( \lambda _{2}-\lambda _{1}\right)
\left( \lambda _{3}-\lambda _{2}\right) }}\left( \left( \lambda
_{2}-\lambda _{1}\right) \left( \lambda _{3}H_{1}-H_{2}\right)
-\allowbreak \left( \lambda _{3}-\lambda _{2}\right) \left(
\lambda _{1}H_{1}-H_{2}\right) \right).\nonumber
\end{eqnarray}
On obtient évidemment des expressions similaires pour les autres
cas. Il serait intéressant de comparer les solutions obtenues sous
forme de séries de Laurent avec les solutions obtenues à l'aide
des fonctions elliptiques de Jacobi (15) ainsi qu'avec celles
obtenues par la méthode des déformations isospectrales (voir
[10,14]).

$b)$ \underline{Cas de Lagrange} : Dans ce cas, on a
$$I_{1}=I_{2},\quad l_{1}=l_{2}=0.$$
Il n'est pas difficile de montrer que dans ce cas aussi,
l'int\'{e}gration s'effectue \`{a} l'aide de fonctions
elliptiques.

$c)$ \underline{Cas de Kowalewski} : Dans ce cas, on a
$$I_{1}=I_{2}=2I_3,\quad l_{3}=0.$$
L'\'{e}tude de ce cas est compliqu\'{e}e. Le syst\`{e}me
diff\'{e}rentiel $\left(13\right),$ s'\'{e}crit explicitement sous
la forme
\begin{eqnarray}
\overset{.}{m}_{1}&=&m_{2}m_{3},\nonumber\\
\overset{.}{m}_{2}&=&-m_{1}m_{3}+2\gamma_{3},\nonumber\\
\overset{.}{m}_{3}&=&-2\gamma_{2},\\
\overset{.}{\gamma}_{1}&=&2m_{3}\gamma _{2}-m_{2}\gamma _{3},\nonumber\\
\overset{.}{\gamma }_{2}&=&m_{1}\gamma
_{3}-2m_{3}\gamma _{1},\nonumber\\
\overset{.}{\gamma }_{3}&=&m_{2}\gamma _{1}-m_{1}\gamma
_{2},\nonumber
\end{eqnarray}
o\`{u}, sans restreindre la g\'{e}n\'{e}ralit\'{e}, nous avons
choisi $l_{2}=0,$ $\mu gl_{1}=1$, $I_{3}=1$ et nous avons
utilis\'{e} la substitution $t\rightarrow 2t$. Ces équations
forment un champ de vecteurs hamiltonien de la forme (9) avec
$n=6$, $m=k=2$, $z=t$, $w= \left( m_{1},m_{2},m_{3}, \gamma_{1},
\gamma_{2}, \gamma_{3}\right)^{\intercal}$,
$$H=\frac{1}{2}\left( m_{1}^{2}+m_{2}^{2}\right) +m_{3}^{2}+2\gamma
_{1},$$ l'hamiltonien et
$$
J=\left(\begin{array}{cccccc}
0&-m_{3}&m_{2}&0&-\gamma_{3}&\gamma _{2}\\
m_{3}&0&-m_{1}&\gamma_{3}&0&-\gamma _{1}\\
-m_{2}&m_{1}&0&-\gamma _{2}&\gamma _{1}&0\\
0&-\gamma_{3}&\gamma_{2}&0&0&0\\
\gamma_{3}&0&-\gamma_{1}&0&0&0\\
-\gamma _{2}&\gamma_{1}&0&0&0&0
\end{array}\right)
$$
En plus des trois int\'{e}grales premi\`{e}res
\begin{eqnarray}
H_{1}&\equiv& H,\nonumber\\
H_{2}&=&m_{1}\gamma _{1}+m_{2}\gamma _{2}+m_{3}\gamma _{3},\nonumber\\
H_{3}&=&\gamma _{1}^{2}+\gamma _{2}^{2}+\gamma _{3}^{2},\nonumber
\end{eqnarray}
le syst\`{e}me ci-dessus admet une quatri\`{e}me int\'{e}grale
premi\`{e}re quartique obtenue par Kowalewski
$$H_{4}=\left( \left( \frac{m_{1}+im_{2}}{2}\right) ^{2}-
\left( \gamma _{1}+i\gamma _{2}\right) \right) \left( \left(
\frac{m_{1}-im_{2}}{2}\right) ^{2}- \left( \gamma _{1}-i\gamma
_{2}\right) \right).$$
Les int\'{e}grales premi\`{e}res $H_{1}$ et
$H_{4}$ sont en involution, tandis que $H_{1}$ et $H_{3}$ sont
triviaux. Soit
$$M_c=\bigcap_{k=1}^{4}\left\{ x:H_{k}\left( x\right) =c_{k}\right\}
,$$ la vari\'{e}t\'{e} affine d\'{e}finie par l'intersection des
quatre constantes du mouvenent. Soit $$( m_{1},m_{2},m_{3},\gamma
_{1},\gamma _{2},\gamma _{3}) \longmapsto (
x_{1},x_{2},m_{3},y_{1},y_{2},\gamma _{3}),$$ une transformation
biunivoque de $M_c$ o\`{u}
$$2x_{1}=m_{1}+im_{2},\qquad y_{1}=x_{1}^{2}-\left( \gamma _{1}+i\gamma _{2}\right)
,$$
$$2x_{2}=m_{1}-im_{2},\qquad y_{2}=x_{2}^{2}-\left( \gamma _{1}-i\gamma _{2}\right)
.$$ Alors, le quotient $M_c/\sigma $ par l'involution
\begin{equation}\label{eqn:euler}
\sigma :M_c\longrightarrow M_c\text{ }\left(
x_{1},x_{2},m_{3},y_{1},y_{2},\gamma _{3}\right) \longmapsto
\left( x_{1},x_{2},-m_{3},y_{1},y_{2},-\gamma _{3}\right) ,
\end{equation}
est une surface $S$ (de Kummer)
\begin{eqnarray}
&&y_{1}y_{2}=c_{4},\nonumber\\
&&y_{1}R(x_{2})+y_{2}R(x_{1})+R_{1}\left( x_{1},x_{2}\right)
+c_{4}(x_{1}-x_{2})^{2}=0,\nonumber
\end{eqnarray}
où
$$R(x)=-x^{4}+c_{1}x^{2}-2c_{2}x+1-c_{4},$$
est un polyn\^{o}me de degr\'{e} $4$ en $x$ et
$$R_{1}(x_{1},x_{2})=-c_{1}x_{1}^{2}x_{2}^{2}-
c_{2}\left( c_{2}-2x_{1}x_{2}\left( x_{1}+x_{2}\right) \right)
+\left( 1-c_{4}\right) \left( c_{1}-\left( x_{1}+x_{2}\right)
^{2}\right),$$ un autre polyn\^{o}me de degr\'{e} $2$\ en
$x_{1},x_{2}.$\ Les points de ramification de $M_c$\ sur $S$\ sont
donn\'{e}s par les points fixes de l'involution $\sigma $ $\left(
17\right) $\ et sont en nombre de $8.$ La surface $S$ est un
rev\^{e}tement double du plan $(x_{1},x_{2}),$ ramifi\'{e} le long
de deux courbes elliptiques se coupant exactement aux $8$ points
fixes de l'involution $\sigma .$ Ces courbes donnent lieu \`{a}
l'\'{e}quation diff\'{e}rentielle d'Euler
$$\frac{dx_{1}}{\sqrt{R(x_{1})}}\pm
\frac{dx_{2}}{\sqrt{R(x_{2})}}=0,$$ \`{a} laquelle se trouvent
li\'{e}s les fameuses variables de Kowalewski
$$s_{1}=\frac{R_{1}\left( x_{1},x_{2}\right)
-\sqrt{R(x_{1})}\sqrt{R(x_{2})}}{\left( x_{1}-x_{2}\right)
^{2}}+\frac{c_{1}}{2},$$
$$s_{2}=\frac{R_{1}\left( x_{1},x_{2}\right)
+\sqrt{R(x_{1})}\sqrt{R(x_{2})}}{\left( x_{1}-x_{2}\right)
^{2}}+\frac{c_{1}}{2},$$ et peuvent \^{e}tre vues comme \'{e}tant
des formules d'addition pour la fonction elliptique de Weierstrass
(voir [8]). En termes des variables $s_{1}$ et $s_{2},$ le
syst\`{e}me $\left(16\right) $ devient
$$\frac{ds_{1}}{\sqrt{P_{5}(s_{1})}}\pm
\frac{ds_{2}}{\sqrt{P_{5}(s_{2})}}=0,$$
$$\frac{s_{1}ds_{1}}{\sqrt{P_{5}(s_{1})}}\pm
\frac{s_{2}ds_{2}}{\sqrt{P_{5}(s_{2})}}=dt,$$ o\`{u} $P_{5}(s)$\
est un polyn\^{o}me de cinqui\`{e}me degr\'{e} et
l'int\'{e}gration s'effectue au moyen des fonctions
hyperelliptiques de genre $2.$

En utilisant le théorème 1, on montre que le syst\`{e}me $\left(
13\right)$ dans le cas de Kowalewski, admet deux familles de
solutions en s\'{e}ries de Laurent m\'{e}romorphes
$$M=\sum_{k=0}^{\infty }M^{\left( k\right) }t^{k-1},\quad \Gamma
=\sum_{k=0}^{\infty }\Gamma ^{\left( k\right)}t^{k-2},$$
d\'{e}pendant de cinq param\`{e}tres libres tels que : les
coefficients $M^{\left( 0\right) }$\ et $\Gamma ^{\left( 0\right)
}$\ satisfont au syst\`{e}me non-lin\'{e}aire
\begin{eqnarray}
&&M^{\left( 0\right) }+\left[ M^{\left( 0\right) },\Lambda
M^{\left( 0\right) }\right] +\left[ \Gamma ^{\left( 0\right)
},L\right] =0,\\
&&2\Gamma ^{\left( 0\right) }+\left[ \Gamma ^{\left( 0\right)
},\Lambda M^{\left( 0\right) }\right] =0,\nonumber
\end{eqnarray}
d\'{e}pendant d'une variable libre $\alpha $\ et d\'{e}finissant
deux droites. Tandis que $M^{\left( k\right) }$\ et $\Gamma
^{\left( k\right) }$\ satisfont aux syst\`{e}mes lin\'{e}aires
$$
\left( L-kI\right)\left(\begin{array}{c}
M^{\left( 1\right) }\\
\Gamma ^{\left( 1\right) }
\end{array}\right)=0,
$$
$$
\left( L-kI\right)\left(\begin{array}{c}
M^{\left( k\right) }\\
\Gamma ^{\left( k\right) }
\end{array}\right)=
\left\{\begin{array}{rl} -\sum_{i=1}^{k-1}\left[ M^{\left(
i\right) },\Lambda M^{\left( k-i\right) }\right] \\
-\sum_{i=1}^{k-1}\left[ \Gamma ^{\left( i\right) },\Lambda
M^{\left( k-i\right) }\right]
\end{array}\right.
,\quad \text{ pour }k\geq 2,
$$
o\`{u} $L$\ est la matrice jacobienne de $\left(18\right) .$\ Ces
syst\`{e}mes fournissent une variable libre \`{a} chacun des
niveaux $k=1,2,3$\ et $4.$ Explicitement, on a\\
$\left( *\right) $ $1^{\grave{e}re}$ famille de solutions en
s\'{e}ries de Laurent m\'{e}romorphes :
\begin{eqnarray}
m_{1}\left( t\right) &=&\frac{\alpha }{t}+i\left( \alpha ^{2}-2\right) \beta +\circ\left( t\right) ,\nonumber\\
m_{2}\left( t\right)& =&\frac{i\alpha }{t}-\alpha ^{2}\beta +\circ\left( t\right) ,\nonumber\\
m_{3}\left( t\right) &=&\frac{i}{t}+\alpha \beta +\circ\left( t\right) ,\nonumber\\
\gamma _{1}\left( t\right)& =&\frac{1}{2t^{2}}+\circ\left( t\right) ,\nonumber\\
\gamma _{2}\left( t\right)& =&\frac{i}{2t^{2}}+\circ\left( t\right) ,\nonumber\\
\gamma _{3}\left( t\right) &=&\frac{\beta }{t}+\circ\left(
t\right) .\nonumber
\end{eqnarray}
$\left( **\right) $ $2^{\acute{e}me}$ famille de solutions en
s\'{e}ries de Laurent m\'{e}romorphes :
\begin{eqnarray}
m_{1}\left( t\right) &=&\frac{\alpha }{t}-i\left( \alpha ^{2}-2\right) \beta +\circ\left( t\right) ,\nonumber\\
m_{2}\left( t\right)& =&-\frac{i\alpha }{t}-\alpha ^{2}\beta +\circ\left( t\right) ,\nonumber\\
m_{3}\left( t\right)& =&-\frac{i}{t}+\alpha \beta +\circ\left( t\right) ,\nonumber\\
\gamma _{1}\left( t\right)& =&\frac{1}{2t^{2}}+\circ\left( t\right) ,\nonumber\\
\gamma _{2}\left( t\right)& =&-\frac{i}{2t^{2}}+\circ\left( t\right) ,\nonumber\\
\gamma _{3}\left( t\right)& =&\frac{\beta }{t}+\circ\left(
t\right) .\nonumber
\end{eqnarray}
Les diviseurs des p\^{o}les des fonctions $M$\ et $\Gamma $\ sont
deux curbes algébriques
\begin{equation}\label{eqn:euler}
\mathcal{D}_{\varepsilon }:\beta ^{4}\left( \alpha ^{2}-1\right)
^{2}-\left( c_{1}\beta ^{2}-2\varepsilon c_{2}\beta -1\right)
\left( \alpha ^{2}-1\right) +c_{4}=0,\text{ }\varepsilon ^{2}=-1,
\end{equation}
irr\'{e}ductibles isomorphes et chacune de genre $3.$\ Ce sont
deux rev\^{e}tements
\begin{equation}\label{eqn:euler}
\mathcal{D}_{\varepsilon }\longrightarrow \mathcal{D}_{\varepsilon
}^{0},\quad \left( \alpha ,u,\beta \right) \longmapsto \left(
u,\beta \right) ,
\end{equation}
doubles ramifi\'{e}s en quatre points de courbes elliptiques :
\begin{equation}\label{eqn:euler}
\mathcal{D}_{\varepsilon }^{0}:u^{2}=\left( c_{1}\beta
^{2}-2\varepsilon c_{2}\beta -1\right) ^{2}-4c_{4}\beta ^{4}.
\end{equation}
Soit $$\mathcal{L}\left( \mathcal{D}_{\varepsilon
=i}+\mathcal{D}_{\varepsilon =-i}\right)= \left\{ f\text{
m\'{e}romorphe sur }:\left( f\right) +\mathcal{D}_{\varepsilon
=i}+\mathcal{D}_{\varepsilon =-i}\geq 0\right\} ,$$ l'espace
vectoriel des fonctions $f$ m\'{e}romorphes telles que: $\left(
f\right) +\mathcal{D}_{\varepsilon =i}+\mathcal{D}_{\varepsilon
=-i}\geq 0.$ , i.e., l'ensemble des fonctions holomorphes en
dehors de $\mathcal{D}_{\varepsilon =i}+\mathcal{D}_{\varepsilon
=-i}$ et ayant au plus des p\^{o}les le long de
$\mathcal{D}_{\varepsilon =i}+\mathcal{D}_{\varepsilon =-i}$. En
utilisant les séries de Laurent obtenues précédemment, on montre
que cet espace est engendr\'{e} par les huit fonctions suivantes:
\begin{eqnarray}
&&f_{0}=1,\text{ }f_{1}=m_{1},\text{ }f_{2}=m_{2},\text{
}f_{3}=m_{3},
\text{ }f_{4}=\gamma _{3},\text{ }f_{5}=f_{1}^{2}+f_{2}^{2},\nonumber\\
&&f_{6}=4f_{1}f_{4}-f_{3}f_{5},\text{ }f_{7}=\left( f_{2}\gamma
_{1}-f_{1}\gamma _{2}\right) f_{3}+2f_{4}\gamma _{2}.\text{ }
\end{eqnarray}
En outre, l'application
$$\left( \mathcal{D}_{\varepsilon =i}+\mathcal{D}_{\varepsilon =-i}\right) \longrightarrow
\mathbb{CP}^{7},p=\left( \alpha ,u,\beta \right) \longmapsto
\underset{t\rightarrow 0}{\lim }t\left(
1,f_{1}(p),...,f_{7}(p)\right)=$$$$ \left( 0,\alpha,\pm i,\pm
i,\beta,\pm i\alpha\beta,\varepsilon(\alpha^2-1)\beta^2,\pm i (\mp
c_2+c_1\beta-2(\alpha^2-1)\beta^3)\right) ,
$$
plonge $\left( \mathcal{D}_{\varepsilon
=i}+\mathcal{D}_{\varepsilon =-i}\right) $ dans $\mathbb{CP}^{7}$\
de telle fa\c{c}on que $\mathcal{D}_{\varepsilon =i}$\ intersecte
$\mathcal{D}_{\varepsilon =-i}$\ transversalement en quatre points
\`{a} l'infini  $\left( \alpha =\pm 1,\text{ }u=\pm \beta
^{2}\sqrt{c_{1}^{2}-4c_{4}},\text{ }\beta =\infty \right)$ et que
le genre g\'{e}om\'{e}trique de $\left( \mathcal{D}_{\varepsilon
=i}+\mathcal{D}_{\varepsilon =-i}\right) $\ est $9.$ Les orbites
du champ de vecteurs en question passant \`{a} travers $\left(
\mathcal{D}_{\varepsilon =i}+\mathcal{D}_{\varepsilon =-i}\right)$
forment une surface lisse $\Sigma $ tout le long de
$\left(\mathcal{D}_{\varepsilon =i}+\mathcal{D}_{\varepsilon
=-i}\right) $ tel que : $\Sigma \backslash
\left(\mathcal{D}_{\varepsilon =i}+\mathcal{D}_{\varepsilon
=-i}\right) \subset M_{c}.$ La vari\'{e}t\'{e}
$\widetilde{M_{c}}=M_{c}\cup \Sigma ,$ est lisse, compacte et
connexe. En outre, le champ de vecteur est r\'{e}gulier le long du
diviseur $\left( \mathcal{D}_{\varepsilon
=i}+\mathcal{D}_{\varepsilon =-i}\right) ,$ transversal en tout
point $\beta \neq 0$ $\left( \text{resp. }\beta \neq \infty
\right) $ et doublement tangent en $\beta =0$ $\left( \text{resp.
}\beta =\infty \right)$. Les champs de vecteurs engendrés par
$H_1$ et $H_4$ se prolongent de fa\c{c}on holomorphe et demeurent
ind\'{e}pendants sur la vari\'{e}t\'{e} $\widetilde{M_{c}}.$ La
vari\'{e}t\'{e} $\widetilde{M_{c}}$ est une surface ab\'{e}lienne
sur laquelle le flot hamiltonien $\left(16\right)$ se
lin\'{e}arise. L'involution $\sigma $ $\left(17\right) $\ sur la
vari\'{e}t\'{e} affine $M_{c}$ se transforme en une involution
$$\sigma _{\varepsilon }:\mathcal{D}_{\varepsilon }\longrightarrow
\mathcal{D}_{\varepsilon },\text{ }\left( \alpha ,u,\beta \right)
\longmapsto \left( -\alpha ,u,\beta \right),$$ sur la surface de
Riemann $\mathcal{D}_{\varepsilon }$ $\left(19\right) $\ et admet
huit points fixes donn\'{e}s par les points de branchements de
$\mathcal{D}_{\varepsilon }$ sur la courbe elliptique
$\mathcal{D}_{\varepsilon }^{0}$ $\left(21\right) .$\ Donc
l'involution en question poss\`{e}de seize points fixes au total,
confirmant ainsi le nombre de points fixes qu'une involution
$z\longmapsto -z$ sur une vari\'{e}t\'{e} ab\'{e}lienne doit en
avoir. Il existe sur la surface ab\'{e}lienne $\widetilde{M_{c}}$
deux diff\'{e}rentielles holomorphes $dt_{1}$ et $dt_{2}$ telles
que :
\begin{eqnarray}
dt_{1\mid _{\mathcal{D}_{\varepsilon }}}&=&\omega
_{1}=\frac{k_{1}\left( \alpha ^{2}-1\right) \beta ^{2}d\beta
}{\alpha u},\nonumber\\
dt_{2\mid _{\mathcal{D}_{\varepsilon }}}&=&\omega
_{2}=\frac{k_{2}d\beta }{\alpha u},\nonumber
\end{eqnarray}
o\`{u} $k_{1},k_{2}\in \mathbb{C}$ et $\omega _{1},\omega _{2}$
sont des diff\'{e}rentielles holomorphes sur
$\mathcal{D}_{\varepsilon }.$ En outre, l'espace des
diff\'{e}rentielles holomorphes sur le diviseur $\left(
\mathcal{D}_{\varepsilon =i}+\mathcal{D}_{\varepsilon =-i}\right)
$ est
$$\left\{ f_{1}^{\left( 0\right) }\omega _{2}\text{ }f_{2}^{\left( 0\right) }
\omega _{2},...,f_{7}^{\left( 0\right) }\omega _{2}\right\} \oplus \left\{ \omega _{1},\omega _{2}\right\}
,$$ o\`{u} $f_{1}^{\left( 0\right) },$\ $f_{2}^{\left( 0\right)
},...$\ $f_{7}^{\left( 0\right) }$\ sont les premiers coefficients
des fonctions $f_{1},f_{2},...,f_{7}$\ $\in \mathcal{L}\left(
\mathcal{D}_{\varepsilon =i}+\mathcal{D}_{\varepsilon =-i}\right)
\left( 66\right) $\ et le plongement de $\left(
\mathcal{D}_{\varepsilon =i}+\mathcal{D}_{\varepsilon =-i}\right)
$\ dans $\mathbb{CP}^{7}$\ est \`{a} deux diff\'{e}rentielles
holomorphes pr\`{e}s le plongement canonique
$$p=\left( \alpha ,u,\beta \right) \in \left( \mathcal{D}_{\varepsilon =i}
+\mathcal{D}_{\varepsilon =-i}\right) \longmapsto \left\{ \omega
_{2},f_{1}^{\left( 0\right) } \omega _{2},f_{2}^{\left( 0\right) }
\omega _{2},...,f_{7}^{\left( 0\right) }\omega _{2}\right\} \in
\mathbb{CP}^{7}.$$ La surface ab\'{e}lienne $\widetilde{M_{c}}$\
est caract\'{e}ris\'{e}e comme \'{e}tant la duale de
vari\'{e}t\'{e} Prym $\left( \mathcal{D}_{\varepsilon
}/\mathcal{D}_{\varepsilon }^{0}\right) $\ du rev\^{e}tement
double $\left(20\right)$. Les solutions sous forme de séries de
Laurent (théorème 1), jouent un rôle crucial dans la preuve de ces
résultats.


\begin{thebibliography}{99}
\bibitem{}ADLER M. AND  VAN MOERBEKE P.: The complex geometry of the Kowalewski-Painlevé
analysis. \emph{Invent. Math.} \textbf{97} (1989) 3-51.
\bibitem{}ARNOLD V.I.: Ordinary differential equations. Springer-Textbook,
3nd edition, 1992.
\bibitem{}CARTAN H.: Théorie élémentaire des fonctions
analytiques d'une ou plusieurs variables complexes. Hermann, 1961.
\bibitem{}FRANCOISE J.P.: Integrability of quasi-homogeneous vector
fields (\emph{preprint}).
\bibitem{}GRIFFITHS P.A. AND HARRIS, J.: Principles of
algebraic geometry. Wiley-Interscience, New York, 1978.
\bibitem{}HAINE L.: Geodesic flow on $SO(4)$\ and Abelian surfaces. \emph{Math.
Ann.} \textbf{263} (1983) 435-472.
\bibitem{}HILLE E.: Ordinary differential equations in the complex
domain. Wiley-Interscience, New-York, 1976.
\bibitem{}LESFARI A.: Abelian surfaces and Kowalewski's top. \emph{Ann. Scient.
\'{E}cole Norm. Sup.}, Paris, sér. 4 \textbf{21} (1988) 193-223.
\bibitem{}LESFARI A.: Eléments d'analyse. Sochepress-Université,
Casablanca, 1991.
\bibitem{}LESFARI A.: Completely integrable systems : Jacobi's heritage. \emph{J.
Geom. Phys} \textbf{31} (1999) 265-286.
\bibitem{}LESFARI A.: Le théorème d'Arnold-Liouville et ses
conséquences. \emph{Elem. Math.} (Issue 1) \textbf{58} (2003)
6-20.
\bibitem{}LESFARI A.: Analyse des singularités de quelques systèmes
intégrables. \emph{C. R. Acad. Sci.} Paris, Ser. \textbf{I 341}
(2005) 85-88.
\bibitem{}LESFARI A.: Abelian varieties, surfaces of general type and integrable
systems. \emph{Beiträge Algebra Geom.}, Vol.48, 1 (2007) 95-114.
\bibitem{}LESFARI A.: Integrable systems and complex geometry. \emph{arXiv :
0706.1579 [math. CV]} (2007) 1-45.
\bibitem{}LESFARI A.: Fonctions et Intégrales elliptiques. \emph{arXiv :
0707.1137 [math. DS]} (2007) 1-46.
\bibitem{}PAINLEVE  P.: Oeuvres: tomes 1,2,3. Edition du C.N.R.S. 1975.
\end{thebibliography}
\end{document}